\documentclass[12pt]{article}
\usepackage{amsmath,amsthm,amssymb,amscd}
\usepackage{latexsym}
\newtheorem{thm}{Theorem}[section]

 \newtheorem{prop}[thm]{Proposition}
 \theoremstyle{definition}
 
 \theoremstyle{remark}

 \numberwithin{equation}{section}

\title
{Exotic $\mathbb{R}^4$'s and positive isotropic curvature}

\author{ Hong Huang}
\date{}
\begin{document}
\maketitle
\begin{abstract}
 We show that no exotic $\mathbb{R}^4$ admits a complete Riemannian metric with uniformly positive isotropic curvature and with bounded geometry.  This is essentially a corollary of the main result in [Hu1], and was stated in [Hu2] without proof. In the process of the proof we also show that the diffeomorphism type of an
infinite connected sum of some connected smooth $n$-manifolds ($n\geq 2$) according to a locally finite graph does not depend on the gluing maps used.


\end{abstract}
\maketitle


\section {Introduction}
In [Hu1] we proved the following result which extends [CZ] to the noncompact case.
\begin{thm} \label{thm 1.1} \ \
 Let $X$ be a
complete, connected, non-compact 4-manifold with uniformly positive isotropic
curvature, with bounded geometry and with no essential incompressible space form. Then $X$ is diffeomorphic to an infinite connected sum of
$\mathbb{S}^4$,  $\mathbb{RP}^4$,  $\mathbb{S}^3\times \mathbb{S}^1$, and /or $\mathbb{S}^3\widetilde{\times} \mathbb{S}^1$.
\end{thm}

Note that here we use the standard smooth structures of $\mathbb{S}^4$,  $\mathbb{RP}^4$,  $\mathbb{S}^3\times \mathbb{S}^1$ and  $\mathbb{S}^3\widetilde{\times} \mathbb{S}^1$.

Now we explain the notion of  infinite connected sum  used here (compare [BBM] and [Hu2]).  Let $G$ be a  countably infinite graph which is connected and locally finite (here we allow an edge to connect a vertex to itself, and  allow more than one edge to connect two vertices (or connect one vertex to itself)),  and let  $\mathcal{X}$ be a class
of connected, smooth $n$-manifolds ($n\geq 2$). We associate an element $X_v\in \mathcal{X}$ to each  vertex $v$ of $G$.  For each edge of $G$, suppose it connects the vertices $v_1$ and $v_2$ (it may be that $v_1=v_2$), we do a connected sum of $X_{v_1}$ and $X_{v_2}$  (as in pp. 102-106 in [BJ]). The result is  a connected, smooth $n$-manifold, which is called an infinite connected sum of members of  $\mathcal{X}$ according to the graph $G$.

The following result is stated in [Hu2] without proof. It is essentially a corollary of Theorem 1.1.

\begin{thm} \label{thm 1.2} \ \
 No exotic $\mathbb{R}^4$ admits a complete Riemannian metric with uniformly positive isotropic curvature and with bounded geometry.
\end{thm}

To prove Theorem 1.2 we also need the fact that the diffeomorphism type of an
infinite connected sum of some connected smooth $n$-manifolds ($n\geq 2$) according to a locally finite graph does not depend on the gluing maps used. This fact is proved in Section 2. Theorem 1.2 itself is proved in Section 3.

\section{Infinite connected sum}

We give more details of the definition of infinite connected sum.  Let $G$ be as in the Introduction. Let $\{v_1, v_2, \cdot\cdot\cdot\}$ be the set of the vertices in $G$. If the connected, smooth $n$-manifold ($n\geq 2$)  $X_{v_i}$ associated to the vertex $v_i$ is orientable, we choose an  orientation of it. For  each pair $(i,j)$ with $i\leq j$, let $\{e_{ij}^k| k=1, 2, \cdot\cdot\cdot, m_{ij}\}$ be the set of edges connecting $v_i$ to $v_j$ (of course, if  there is no edge
 connecting $v_i$ to $v_j$, $m_{ij}=0$, and in this case this set is empty).  To each edge $e_{ij}^k$ we associate a pair of smooth
  embeddings $f_{e_{ij}^k}: \mathbb{R}^n\rightarrow X_{v_i}$ and $g_{e_{ij}^k}: \mathbb{R}^n\rightarrow X_{v_j}$ from the standard $\mathbb{R}^n$.
We fix an orientation of the standard $\mathbb{R}^n$. If $X_{v_i}$ is oriented, we let $f_{e_{ij}^k}$ ($j\geq i$, $k=1,2,\cdot\cdot\cdot, m_{ij}$) be orientation-preserving, and let $g_{e_{pi}^l}$ ($p\leq i$, $l=1,2, \cdot\cdot\cdot, m_{pi})$  be orientation-reversing.
We assume that the images of all these embeddings are disjoint from each other.

Let $D^n$ be the closed unit ball with center the origin in the standard $\mathbb{R}^n$.

For each $i$, let
\begin{equation*}
Y_i:=X_{v_i}\setminus (\cup_{j\geq i}\cup_{k=1}^{m_{ij}}f_{e_{ij}^k}(\frac{1}{3}D^n)\cup\cup_{p\leq i}\cup_{l=1}^{m_{pi}}g_{e_{pi}^l}(\frac{1}{3}D^n)),
\end{equation*}
and let $Y$ be the infinite  disjoint union $\bigsqcup Y_i$.

We define an equivalence relation $\sim$ in  $Y$
 by setting
\begin{equation*}
f_{e_{ij}^k}(tu)\sim g_{e_{ij}^k}((1-t)u)
\end{equation*}
 for all $(i,j)$ with $i\leq j$, $k=1,2,\cdot\cdot\cdot, m_{ij}$, $t \in (\frac{1}{3},\frac{2}{3})$ and $u \in \mathbb{S}^{n-1}$. Let $X$ be the quotient space $Y/\sim$.  We call $X$ the connected sum of $X_{v_i}$ according to the graph $G$ via $\{f_{e_{ij}^k}, g_{e_{ij}^k}\}$, and denote it by  $\sharp_G X_{v_i}(f_{e_{ij}^k}, g_{e_{ij}^k})$.

Let $M^n$ be a connected $n$-manifold, and $\varphi=\sqcup_{i=1}^k  \varphi_i: \sqcup_{i=1}^k D^n\rightarrow M^n$ and $\tilde{\varphi}=\sqcup_{i=1}^k  \tilde{\varphi}_i: \sqcup_{i=1}^k D^n\rightarrow M^n$
be two embeddings from the disjoint union of $k$ copies of the standard $n$-disk.
As in [BJ, Definition (10.1)] we say $\varphi$ and $\tilde{\varphi}$ are compatibly oriented if either $M^n$ is not orientable, or, for each $i$ ($1\leq i \leq k$), $\varphi_i$ and $\tilde{\varphi}_i$ are both orientation preserving or both orientation reversing (relative to fixed orientations of $D^n$ and $M^n$).

The following result is well-known, see for example Theorem 3.2 in Chapter 8 of [H].
\begin{prop} \label{prop 2.1} \ \
Let  $\varphi=\sqcup_{i=1}^k  \varphi_i: \sqcup_{i=1}^k D^n\rightarrow M^n$ and $\tilde{\varphi}=\sqcup_{i=1}^k  \tilde{\varphi}_i: \sqcup_{i=1}^k D^n\rightarrow M^n$
be two (smooth) embeddings from the disjoint union of $k$ ($< \infty$) copies of the standard $n$-disk to a connected, smooth $n$-manifold $M^n$ ($n\geq 2$). Suppose that $\varphi$ and $\tilde{\varphi}$ are compatibly oriented. Then there is a diffeotopy $H$ of $M^n$, which is fixed outside of a compact subset of $M^n$ such that $H(\cdot, 1)\circ \varphi= \tilde{\varphi}$.
\end{prop}

\noindent (For definition of diffeotopy (or ambient isotopy), see  [BJ, Definition (9.3)] and p.178 of [H].)

{\bf Proof}  We follow closely the proof of Theorem 3.2 in Chapter 8 of [H]. We do induction on $k$. The $k=1$ case is due to Cerf and Palais (for expositions see Chapters 9 and 10 in [BJ], Chapter III in [K] and Theorem 3.1 in Chapter 8 of [H]). Suppose the result is true for $k=j$. Now we consider the case $k=j+1$. By assumption there exists a diffeotopy $\widetilde{H}$ of $M^n$, which is fixed outside of a compact subset of $M^n$ such that  $\widetilde{H}(\cdot, 1)\circ \varphi|\sqcup_{i=1}^j D^n = \tilde{\varphi}|\sqcup_{i=1}^j D^n $. (In particular, it follows that $\widetilde{H}(\cdot,1)(\varphi_{j+1}(D^n)) \subset M^n \setminus \cup_{i=1}^j \tilde{\varphi}_i(D^n)$.)
Since $n\geq 2$, $M^n \setminus \cup_{i=1}^j \tilde{\varphi}_i(D^n) $ is connected.  We apply  the $k=1$ case to the two embeddings
\begin{equation*}
\widetilde{H}(\cdot,1)\circ \varphi_{j+1}, \tilde{\varphi}_{j+1}: D^n \rightarrow M^n \setminus \cup_{i=1}^j \tilde{\varphi}_i(D^n),
\end{equation*}
and get  a diffeotopy $\hat{H}$ of $M^n \setminus \cup_{i=1}^j \tilde{\varphi}_i(D^n)$ which is fixed outside of a compact subset of  $M^n \setminus \cup_{i=1}^j \tilde{\varphi}_i(D^n)$ such that  $\hat{H}(\cdot, 1)\circ \widetilde{H}(\cdot,1)\circ \varphi_{j+1} =\tilde{\varphi}_{j+1}$. Clearly $\hat{H}$ extends to a diffeotopy of $M^n$ which leaves $\cup_{i=1}^j \tilde{\varphi}_i(D^n)$ fixed. Then $H_t:=\hat{H}_t\circ \widetilde{H}_t$ is the desired diffeotopy. \hfill{$\Box$}

\begin{thm} \label{thm 2.2} \ \ The infinite connected sum
   $\sharp_G X_{v_i}(f_{e_{ij}^k}, g_{e_{ij}^k})$ is a connected, smooth manifold, and oriented if all $X_{v_i}$ are oriented. Its diffeomorphism type (oriented if relevant) does not depend on the choice of embeddings $f_{e_{ij}^k}$ and $g_{e_{ij}^k}$.
\end{thm}
{\bf Proof} The first claim can be shown  as in pp. 103-104 in [BJ] and pp. 90-91 in [K]. Now we show the second claim. Suppose that to each edge $e_{ij}^k$ we associate another pair of smooth
  embeddings $\tilde{f}_{e_{ij}^k}: \mathbb{R}^n\rightarrow X_{v_i}$ and $\tilde{g}_{e_{ij}^k}: \mathbb{R}^n\rightarrow X_{v_j}$ from the standard $\mathbb{R}^n$. If $X_{v_i}$ is oriented, we let $\tilde{f}_{e_{ij}^k}$ ($j\geq i$, $k=1,2,\cdot\cdot\cdot, m_{ij}$) be orientation-preserving, and let $\tilde{g}_{e_{pi}^l}$ ($p\leq i$, $l=1,2, \cdot\cdot\cdot, m_{pi})$  be orientation-reversing. We assume that the images of all these embeddings $\tilde{f}_{e_{ij}^k}$ and $\tilde{g}_{e_{pi}^l}$  are disjoint from each other. We define $\widetilde{Y}_i$ and $\widetilde{Y}$ as before using $\tilde{f}_{e_{ij}^k}$ and  $\tilde{g}_{e_{ij}^k}$. We also introduce an equivalence relation in $\widetilde{Y}$ as before using $\tilde{f}_{e_{ij}^k}$ and  $\tilde{g}_{e_{ij}^k}$, and still denote it by $\sim$. Finally we let  $\sharp_G X_{v_i}(\tilde{f}_{e_{ij}^k}, \tilde{g}_{e_{ij}^k}):= \widetilde{Y}/\sim$.  For each $i$ let
\begin{equation*}
\varphi_i:=\sqcup_{j\geq i} \sqcup_{k=1}^{m_{ij}}  f_{e_{ij}^k} \sqcup  \sqcup_{p\leq i} \sqcup_{l=1}^{m_{pi}} g_{e_{pi}^l}: \sqcup \mathbb{R}^n\rightarrow X_{v_i}
\end{equation*}
and
\begin{equation*}
\tilde{\varphi}_i:=\sqcup_{j\geq i} \sqcup_{k=1}^{m_{ij}}  \tilde{f}_{e_{ij}^k} \sqcup  \sqcup_{p\leq i} \sqcup_{l=1}^{m_{pi}} \tilde{g}_{e_{pi}^l}: \sqcup \mathbb{R}^n\rightarrow X_{v_i}.
\end{equation*}

 Since $G$ is locally finite, for each $i$, the above $\sqcup \mathbb{R}^n$ is a finite disjoint union; we consider the finite disjoint union $\sqcup D^n$ contained in it. For each $i$, we can apply Proposition 2.1 to $\varphi_i|\sqcup D^n$ and $\tilde{\varphi}_i|\sqcup D^n$, and get a  diffeotopy $H_i(\cdot,t)$ of $X_{v_i}$, which is fixed outside a compact subset of $X_{v_i}$, such that
\begin{equation}
\tilde{\varphi}_i|\sqcup D^n=H_i(\cdot,1)\circ \varphi_i|\sqcup D^n.
\end{equation}
By equation (2.1)  we can define a map $F: Y=\sqcup Y_i \rightarrow \widetilde{Y}=\sqcup \widetilde{Y}_i$ via
\begin{equation*}
F(x)=H_i(x,1)  \hspace{4mm} when   \hspace{4mm} x \in Y_i \hspace{4mm}  for \hspace{4mm} some \hspace{4mm} i.
 \end{equation*}
 Clearly $F$ is a diffeomorphism. Note that by equation (2.1) again $F$ is compatible with  the equivalence relations in $Y$ and in $\widetilde{Y}$.  So $F$ induces a diffeomorphism

\begin{equation*}
 \overline{F}: \sharp_G X_{v_i}(f_{e_{ij}^k}, g_{e_{ij}^k})\rightarrow \sharp_G X_{v_i}(\tilde{f}_{e_{ij}^k}, \tilde{g}_{e_{ij}^k}).
\end{equation*}
\hfill{$\Box$}

\vspace*{0.4cm}

\noindent {\bf Remark} In general, if each $X_{v_i}$ is orientable, the diffeomorphism type of the infinite connected sum $\sharp_G X_{v_i}$ may depend on  the choice of the orientations of $X_{v_i}$. But it is easy to see that if each $X_{v_i}$ is orientable and admits an orientation-reversing diffeomorphism, then   the (unoriented) diffeomorphism type of $\sharp_G X_{v_i}$ does not depend on the choice of the orientations of $X_{v_i}$.

\section{ Proof of Theorem 1.2}

Let $X$ be a smooth 4-manifold which is homeomorphic to the standard $\mathbb{R}^4$. Assume that $X$ admits a complete Riemannian metric with uniformly positive isotropic curvature and with bounded geometry.  Clearly $X$ contains no essential incompressible space form. By  Theorem 1.1,
$X$ is diffeomorphic to an infinite connected sum of
$\mathbb{S}^4$,  $\mathbb{RP}^4$,  $\mathbb{S}^3\times \mathbb{S}^1$, and /or $\mathbb{S}^3\widetilde{\times} \mathbb{S}^1$ according to a locally finite graph $G$. Since the fundamental group of $X$ is trivial, the graph $G$ must be a tree, and the smooth manifold $X_v$ associated to any vertex $v$ in $G$ must be diffeomorphic to the standard $\mathbb{S}^4$.

We know that any topological manifold homeomorphic to $\mathbb{R}^4$ has exactly  one (topological) end (for definition see for example, [DK]). It follows that the tree $G$ has exactly  one topological end also. But for a locally finite graph, there is a natural
bijection between its topological ends and its graph-theoretical ends, cf. [DK] and the references therein.  So the tree $G$ has only one graph-theoretical end. (There should be a more direct argument for this fact.) Now we choose a ray $\gamma$ in $G$, which is essentially unique.
Let $w_0, w_1, w_2, \cdot\cdot\cdot$ be the set of vertices along the ray $\gamma$.  For each $i$, there are only finite vertices which can be connected to $w_i$ via a sequence of edges not contained in the ray $\gamma$.
For each $i$, we do connected sum of all $\mathbb{S}^4$'s associated to these finite vertices (including $w_i$), the result is diffeomorphic to the $\mathbb{S}^4$ associated to the vertex $w_i$ via a diffeomorphism not affecting the part of this $\mathbb{S}^4$ where its connected sum with the two $\mathbb{S}^4$'s associated to $w_{i-1}$ and $w_{i+1}$ occurs. Then we see that $X$ is diffeomorphic to an infinite connected sum of $\mathbb{S}^4$'s according to the ray $[0, +\infty)$ with a vertex $w_i$ at $i$ ($i=0, 1, 2, \cdot\cdot\cdot$) using some gluing maps. 

We know that the infinite connected sum of $\mathbb{S}^4$'s (all with the standard orientation) according to the ray $[0, +\infty)$  using certain special gluing maps actually produces the standard $\mathbb{R}^4$. (Represent the standard $\mathbb{R}^4$ as the union of the unit 4-ball and the closed subspaces $A_i$ bounded by the two 3-spheres with radius $i$ and $i+1$ (and each with center the origin), $i=1, 2, \cdot\cdot\cdot$.  Note that each $A_i$ may be seen as $\mathbb{S}^4$ with two open 4-balls removed.)  We also know that $\mathbb{S}^4$ admits an orientation-reversing diffeomorphism. So by Theorem 2.2 and the Remark following it, $X$ is diffeomorphic to the standard $\mathbb{R}^4$.   \hfill{$\Box$}

\vspace*{0.4cm}

\noindent {\bf Remark} It is interesting to see whether or not the condition `with bounded geometry' in Theorem 1.2 can be removed.

\vspace*{0.4cm}

\noindent {\bf Acknowledgements} {\hspace*{4mm}}   This note was conceived during the workshop on the geometry of submanifolds  and curvature flows, April 25-29, 2016 in Zhejiang University. I would like to thank all the organizers and  participants of this workshop, in particular, Prof. Zizhou Tang for his kind invitation, Prof. Hongwei Xu for his hospitality, Prof. Bing-Long Chen for his questions on my talk  in this workshop, and Prof. Weiping Zhang for his encouragements.


\hspace *{0.4cm}

\bibliographystyle{amsplain}

\noindent {\bf References}

\hspace *{0.1cm}

\vspace*{0.4cm}

[BBM] L. Bessi$\grave{e}$res, G. Besson and S. Maillot, Ricci flow on
 open 3-manifolds and positive scalar curvature,  Geometry and Topology  15 (2011), 927-975.

[BJ]  Th. Br$\ddot{o}$cker and K. J$\ddot{a}$nich,  Introduction to differential topology, Cambridge University Press 1982.

[CZ] B.-L. Chen, X.-P. Zhu, Ricci flow with surgery
on  four-manifolds with positive isotropic curvature, J. Diff. Geom.
74 (2006), 177-264.

[DK] R. Diestel, D. K$\ddot{u}$hn,  Graph-theoretical versus topological ends of graphs, J. Combin. Theory Ser. B 87 (2003), no.1, 197-206.

[H] M. Hirsch, Differential topology, GTM 33, Springer-Verlag 1994.

[Hu1] H. Huang, Ricci flow on open 4-manifolds with positive isotropic
curvature, J. Geom. Anal. 23 (2013), no.3, 1213-1235.

[Hu2] H. Huang, Four-orbifolds with positive isotropic curvature, Comm. Anal. Geom. 23 (2015), no.5, 951-991.

[K] A. Kosinski, Differential manifolds, Academic Press 1993.

\vspace *{0.4cm}

School of Mathematical Sciences, Beijing Normal University,

Laboratory of Mathematics and Complex Systems, Ministry of Education,

Beijing 100875, P.R. China

 E-mail address: hhuang@bnu.edu.cn

\end{document}